\documentclass[final,nomarks]{dmtcs-episciences}

\usepackage[utf8]{inputenc}
\usepackage{subfigure}

\usepackage{amsmath,amsthm}
\usepackage{amssymb}

\usepackage{color}
\usepackage{xspace}
\usepackage{tikz}
\usetikzlibrary{calc}

\definecolor{block}{RGB}{255,255,255}
\newenvironment{blockmatrix}{%
\vcenter\bgroup\hbox\bgroup
  \tikzpicture[
    x=1.5\baselineskip,
    y=1.5\baselineskip,
  ]%
}{%
  \endtikzpicture
  \egroup
  \egroup
}

\newcommand*{\block}[1][block]{%
  \blockaux{#1}%
}
\def\blockaux#1(#2,#3)#4(#5,#6){%
  \draw[fill={#1}]
  let \p1=(#2,#3),
      \p2=(#5,#6),
      \p3=(#2+#5,#3+#6),
      \p4=(#2+#5/2,#3+#6/2)
  in
    (\p1) rectangle (\p3)
    (\p4) node {$#4$}
  ;%
}

\def\sk{ \mskip3mu }
\def\v{\vert}
\def\a{\ensuremath{\mathcal A}\xspace}

\def\ss{start-small\xspace}
\def\st{\textrm{St}}
\def\si{\sigma}
\def\rl{right-to-left max\xspace}
\def\av{$\{1243,\,2134\}$-avoider\xspace}
\def\avs{$\{1243,\,2134\}$-avoiders\xspace}
\def\gf{generating function\xspace}
\def\gfs{generating functions\xspace}

\def\mbf#1{\mathchoice{\hbox{\boldmath $\displaystyle #1$}}
        {\hbox{\boldmath $\textstyle #1$}}
        {\hbox{\boldmath $\scriptstyle #1$}}
        {\hbox{\boldmath $\scriptscriptstyle #1$}}} 

\newtheorem{theorem}{Theorem}

\newtheorem{lemma}[theorem]{Lemma}
\newtheorem{prop}[theorem]{Proposition}

\author{David Callan}
\title{The number of $\{1243,\,2134\}$-avoiding permutations}
\affiliation{Department of Statistics,  University of Wisconsin-Madison, USA}
\received{2019-03-15}
\revised{2022-07-25, 2023-01-18}
\accepted{2023-03-06}
\begin{document}
\publicationdetails{22}{2023}{2}{13}{5287}
\maketitle
\begin{abstract}
~

We show that the counting sequence for permutations avoiding both of the (classical) patterns 1243 and 2134 
has the algebraic generating function supplied by Vaclav Kotesovec for sequence A164651 in The On-Line Encyclopedia of Integer Sequences.
\end{abstract}
%

\section{Introduction} 

Several authors have developed methods to count permutations avoiding a given set of patterns; see, for example, the references in the Wikipedia entry \cite{wikiPattern}.
In particular, Vatter \cite{schemes2008} and Pudwell \cite{schemes2010} have devised enumeration schemes for automated counting. When successful, an automated method produces an enumeration scheme that yields the initial terms of the counting sequence (perhaps 20 or more) and sometimes these terms appear to have an algebraic \gf that does not follow readily from the enumeration scheme. Here, we treat one such case. The counting sequence for permutations that avoid both of the (classical)
patterns 1243 and 2134 begins $1,\,2,\, 6,\, 22,\, 87,\, 354,\, \dots$, sequence 
\htmladdnormallink{A164651}{http://oeis.org/A164651}
in The On-Line Encyclopedia of Integer Sequences \cite{oeis}. In a comment on this sequence
dated Oct 24 2012, Vaclav Kotesovec observed that the \gf 
\[
\frac{3x^2-9x+2+x(1-x)\sqrt{1-4x}}{2(x-1)(x^2+4x-1)}
\]
fits the known terms of the sequence.

We show that $\{1243,\,2134\}$-avoiders do indeed have this \gf. Defining a \ \textit{start-small} 
permutation to be one that does not start with its largest entry,
the proof rests on a bijection $\phi$ from start-small \avs of length $n$ to lists of \ss 123-avoiders  
whose total length is $n-1$ plus the length of the list. Note that a start-small permutation has length at least 2.   
The Catalan numbers count 123-avoiders and the  
combinatorial interpretation of the Invert transform
permits counting these lists of \ss 123-avoiders. The \gf for start-small \avs thus 
obtained readily yields the \gf for unrestricted \avs. 

In Section 2, we show that $\phi$ is given by iteration of a more basic bijection $\psi$, which is presented in Section 3. Finally, Section 4 gives the bookkeeping details to obtain the desired \gf from~$\phi$.

\section{The bijection $\mbf{\phi}$}
A \textit{mid-123} entry in a permutation is an entry that serves as the ``2'' in a 123 pattern. A \ \textit{key} mid-123 entry is 
a mid-123 entry $b$ whose immediate predecessor is either less than $b$ or a right-to-left maximum (max for short). For example, the mid-123 entries in
1\,6\,3\,2\,4\,5 are 3,\,2,\,4 but only 3 and 4 are key.

\begin{lemma} \label{key}
A permutation with no key mid-$123$ entries has no mid-$123$ entries at all and so is a $123$-avoider.
\end{lemma}
\begin{proof} 
We suppose that $\pi$ has no key mid-123 entry and $\pi_j$ is the 
first mid-123 entry, and deduce a contradiction. Say $\pi_i\pi_j$ begins a 123 pattern. Since $\pi_j$ is not key, we have (i) $\pi_{j-1}>\pi_j$ 
and (ii) $ \pi_{j-1}$ is not a right-to-left max. 
Now (i) implies $i\ne j-1$ and (ii) implies $\pi_{\ell}>\pi_{j-1}$ 
for some $\ell>j-1$. But then $\pi_i \pi_{j-1}\pi_{\ell}$ is a 123, contradicting the minimality of $j$.
\end{proof} 

\begin{lemma} \label{uniquec}
Suppose $b$ is the last mid-$123$ entry in a \av. Then there is only one entry that occurs after $b$ in $\pi$ and is larger than $b$.
\end{lemma}
\begin{proof}  Suppose $abc$ and $abc'$ are 123 patterns with $c\ne c'$, say $c<c'$. If $c$ precedes $c'$ in the permutation, then $c$ is a mid-123, violating the hypothesis on $b$. If $c$ follows $c'$, then $abc'c$ is an instance of a forbidden 1243 pattern. 
\end{proof} 

Let $\a_{n}$ denote the set of \ss \avs on $[n]$, and $\a_{n,k}$ the subset with $k$ key mid-123 entries. 
To produce the promised bijection $\phi$ from $\a_{n}$ to lists of \ss 123-avoiders whose total length 
is $n-1$ plus the length of the list, it suffices to exhibit, for $0\le k \le n-2$, a bijection 
\begin{eqnarray*}
\a_{n,k} \rightarrow \big\{(\si_{1}, \dots, \si_{k+1}): &\textrm{each $\si_{i}$ is a \ss 123-avoider,} \\[-1mm] 
 & \textrm{lengths of the $\si_{i}$'s sum to $n+k$}\big\}.
\end{eqnarray*}
For $k=0$, naturally the bijection is $\pi   \rightarrow (\pi)$, a singleton list, because, by Lemma \ref{key}, 
$\pi$ is already 123-avoiding.
For $1 \le k \le n-2$ and $k+1 \le j \le n-1$ let 
\[
\a_{n,k,j}=\{\pi \in \a_{n,k}:\ \textrm{the last mid-123 entry of $\pi$ is in position $j$}\}.
\] 
In the next section, for $1\le k < j < n$, we present a bijection $\psi$ between $\a_{n,k,j}$ and $ \a_{j,\,k-1}  \times  \a_{n+1-j,\,0}$. Iteration of $\psi$ then gives the required bijection $\phi$.

\section{The bijection $\psi$}
To \textit{standardize}  a list $\tau$ of distinct integers means to replace the smallest by 1, the next smallest by 2, and so on, and the result is denoted \st($\tau$).
We define $\psi$ as follows: given $\pi \in \a_{n,k,j}$, we need to reversibly produce $\si_{1} \in \a_{j,\,k-1}$ and $\si_{2} \in \a_{n+1-j,\,0}$.
Write $\pi$ as $\tau_{1}b\tau_{2}$ where $b$ is the last mid-123 entry in $\pi$, and let $abc$ be the 123 pattern in $\pi$ with smallest $a$ ($c$ is uniquely determined by Lemma \ref{uniquec}). Concatenate $a$ and $\tau_{2}$ and standardize to get the desired $\si_{2}=\st(a\tau_{2})$ with no 123's. 

Concatenate  $\tau_{1}$ and $c$ to get a \av $\rho$---a candidate (after standardization) for $\si_{1}$. 
This $\rho$ may need further processing because of two glitches: $\rho$ may still 
have $k$ key mid-123's instead of the required $k-1$ and $\rho$ cannot end with its 
smallest entry (which must be possible in $\si_{1}$). But these glitches cancel out. 

If $b$ is a key mid-123 in $\pi$, then $\rho$ has $k-1$ key mid-123's because $b$ 
has been lost. In this case, just standardize $\rho$ to get $\si_{1}$. 
If $b$ is not key, the last entry of $\tau_{1}$ exceeds $b$ but is not a \rl of $\pi$. Delete from $\rho$ the longest right factor (terminal 
string) of $\tau_{1}$ that is decreasing but does not start with a \rl of $\pi$, 
equivalently, does not start with an entry $>c$; say $t \ge 1$ entries are deleted. Add $t$ to each remaining 
entry of $\rho$, append the entries $t,t-1, \dots, 1$ and standardize to obtain $\si_{1}$. 

Before proving the invertibility of $\psi$, we give 3 examples according as 
\begin{enumerate}
 \item $\pi_j=b$ is key and $\pi_{j-1}<b$,
 \item $\pi_j$ is key and $\pi_{j-1}$ is a \rl,
 \item $\pi_j$ is not key.
\end{enumerate}

\begin{figure}[h!]
\begin{center}
\includegraphics[angle=0, scale =0.8]{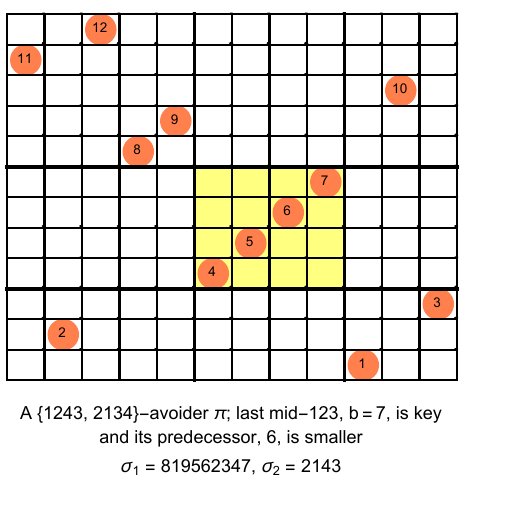}
\includegraphics[angle=0, scale =0.8]{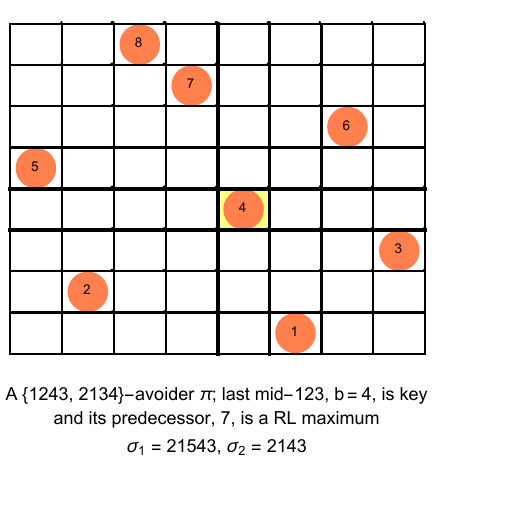}

\includegraphics[angle=0, scale = 0.8]{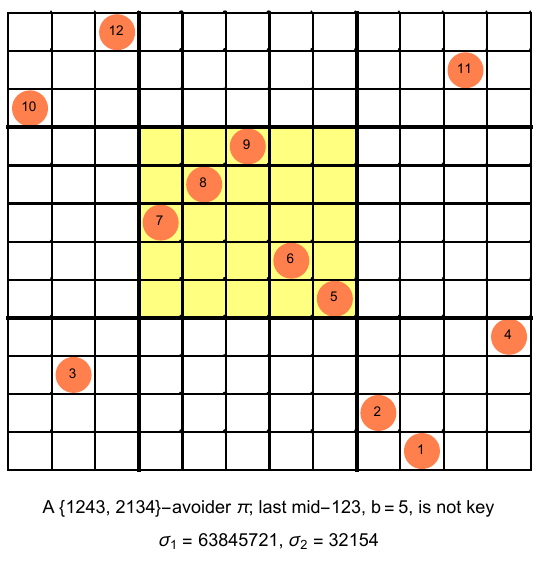}
\caption{Three examples of the action of $\psi$. \label{fig1}} 
\end{center}
\end{figure}

These examples are illustrated in Figure \ref{fig1} below (the yellow box in the figures is explained in Figure \ref{fig2} further below). 
\begin{enumerate}
 \item $\pi= 11\mskip6mu 2\mskip6mu 12\mskip6mu 8 \sk  9 \sk  4 \sk  5 \sk  6 \sk  7 \sk  1\mskip5mu 10\mskip6mu 3$ decomposes as $\tau_1 \mskip6mu 7\mskip6mu  \tau_2$ with $a=2,\mskip6mu b=7,\mskip6mu c=10$. So $\si_1=\st(11\mskip6mu 2\mskip6mu 12\mskip6mu 8 \sk  9 \sk  4 \sk  5 \sk  6\mskip6mu 10)= 8 \sk  1 \sk  9 \sk  5 \sk  6 \sk  2 \sk  3 \sk  4 \sk  7$ and $\si_2=\st(2 \mskip4mu   1\mskip6mu 10\mskip6mu 3)=2 \sk  1 \sk  4 \sk  3$.
 \item $\pi=5 \sk  2 \sk  8 \sk  7 \sk  4 \sk  1 \sk  6 \sk  3$ decomposes as $\tau_1 \mskip6mu 4\mskip6mu  \tau_2$ with $a=2,\mskip6mu b=4,\mskip6mu c=6$  So $\si_1=\st(5 \sk  2 \sk  8 \sk  7 \sk  6)=2\sk 1\sk 5\sk 4\sk 3$ and $\si_2=\st(2\sk 1\sk 6\sk 3)=2\sk 1\sk 4\sk 3$.
 \item $\pi=10\mskip6mu 3\mskip6mu  12\mskip6mu 7 \sk  8 \sk  9 \sk  6 \sk  5 \sk  2 \sk  1\mskip6mu 11\mskip6mu  4$ decomposes as $\tau_1 \mskip6mu 5\mskip6mu  \tau_2$ with $a=3,\mskip6mu b=5,\mskip6mu c=11$. Here, $b$ is not key and $t=2$, so $\si_1=\st((2+(10\mskip6mu  3\mskip6mu  12\mskip6mu  7 \sk  8 \sk  11)) \sk  2 \sk  1) = 6\sk 3\sk 8\sk 4\sk 5\sk 7\sk 2\sk 1$ and $\si_2=\st(3 \sk  2 \sk  1\mskip5mu 11\mskip6mu 4)=3\sk 2\sk  1\sk  5\sk 4$.
\end{enumerate}

Next, we investigate the structure of \ss $\{1243,\,2134\}$-avoiders $\pi$  to obtain a more explicit description of $\psi$. Recall $\pi_j=b$ is the last mid-123 entry in $\pi$, and $abc$ is the 123 pattern in $\pi$ with smallest $a$ and $c$ is uniquely determined by Lemma \ref{uniquec}. 
Recall also that $b$ is key if $\pi_{j-1}<b$ or  $\pi_{j-1}$ is a \rl.
Define $i$ and $k$ by $\pi_i=a$ and $\pi_k=c$.

\begin{lemma}\label{bKey} Suppose $b$ is key. Take $\ell$ minimal in $[i+1,j]$ such that $\pi_\ell<\pi_{\ell+1}\dots<\pi_j\ (\ell=j$ if $\pi_{j-1}$ is a \rl\!$)$. Then 
\begin{enumerate}[(i)]
 \item $\pi_\ell,\dots,\pi_j$ are consecutive integers, and 
 \item For $h<\ell, \ \pi_h<b$ implies $h=i$.
\end{enumerate}
\end{lemma}
\begin{proof} $(i)$ If not, say $\pi_h<x<\pi_{h+1}$ with $\ell\le h <j$ and consider the position of $x$ in $\pi$. If $x$ occurs after $c$, then $a\pi_hcx$ is a forbidden 1243. If $x$ occurs between $b$ and $c$, then $axc$ is a 123 contradicting the assumption that $b$ is the last mid-123. Otherwise, since $\pi_h,\pi_{h+1},\dots,\pi_j=b$ are contiguous entries in $\pi$, $x$ must occur before $\pi_h$ and $x\pi_h\pi_{h+1}c$ is a forbidden 2134. 

$(ii)$ If not, by $(i)$ there exists $h<\ell$ with $\pi_h<\pi_\ell$ and $h\ne i$. Now 
$\pi_{\ell-1}=a$ or $\pi_{\ell-1}>b$ by the definition of $\ell$, and $\pi_h>a$ by the definition of $a$. In case $\pi_{\ell-1}=a$, $\pi_h a bc$ is a 2134, and in case $\pi_{\ell-1}>b$, $\pi_h abc$ is a 2134 if $\pi_h$ occurs before $a$ and $a\pi_h \pi_{\ell-1}b$ is a 1243 if $\pi_h$ occurs after $a$. 
\end{proof}
\begin{lemma} \label{bNotKey} 
Suppose $b$ is not key. Take $r \ge 1$ maximal such that $\pi_{j-r}$ is not a \rl and 
$\pi_{j-r}>\pi_{j-r+1}>\dots>\pi_j$. Then 
\begin{enumerate}[(i)]
 \item $\pi_j,\pi_{j-1},\dots,\pi_{j-r+1}$ are consecutive integers, and 
 \item the integers in the interval $\pi_{j-r+1}+1,\dots,\pi_{j-r}-1$ $($if any$)$ all occur immediately before $\pi_{j-r}$ and in increasing order.
\end{enumerate}
\end{lemma}
\begin{proof}
Take $u\ge 0$ maximal such that $\pi_j,\pi_{j-1},\dots,\pi_{j-u}$ is an increasing sequence 
of consecutive integers $b,b+1,\dots,b+u$. Thus $\pi_{j-(u+1)}\ne b+u+1$. 
If $u\ge 1$ and either $\pi_{j-(u+1)}<\pi_{j-u}$ or $\pi_{j-(u+1)}$ is  
a  right-to-left max, then the result holds with $r=u$  and no integers in the interval $[\pi_{j-r+1}+1,\pi_{j-r}-1]=[b+u,b+u-1]$ of part $(ii)$.
Otherwise, using the fact $b$ is not key when $u=0$, we have $\pi_{j-(u+1)}>b+u+1$, say 
$\pi_{j-(u+1)}=b+u+s$ with $s\ge 2$, and $b+u+s<c$ since $\pi_{j-(u+1)}$ is not a  right-to-left max. Also, $a$ occurs before $b+u+1$ for else $(b+u+1)\,a\,(b+u+s)\,c$ is a 2134. 
The integers $b+u+1,\dots,b+u+s-1$ all occur before $\pi_{j-(u+1)}=b+u+s$ by Lemma \ref{uniquec}, 
and in increasing order for else $(b+u+s)\,c$ is the ``34'' of a 2134 pattern. 
Furthermore, $b+u+1,\dots,b+u+s-1$ all lie directly to the left of $b+u+s$ 
for else some $x \notin [b,b+u+s]$ lies between $b+u+1$ and $b+u+s$ in $\pi$; if $x<b,$ then $(b+u+1)\,x\,(b+u+s)\,c$ is a 2134 while if $x>b+u+s$, then $a\,(b+u+1)\,x\,(b+u+s)$ is a 1243. The result follows with $r=u+1$ and a nonempty interval in part $(ii)$.
\end{proof}
\begin{figure}[ht]
\begin{center}
$\begin{blockmatrix}
    \block(0,0)A(1,1)
    \block[yellow](1,1)B(1,1)
    \block(2,2)C(1,1)
    \block(0,1)(1,1)
    \block(0,2)X(1,1)
    \block(1,2)(1,1)
    \block(1,0)(1,1)
    \block(2,1)(1,1)
    \block(2,0)Z(1,1)
\end{blockmatrix}$
\end{center}
\caption{Matrix form of a $\{1243, 2134\}$-avoider. \label{fig2}}
\end{figure}

It follows from Lemmas \ref{bKey} and \ref{bNotKey} that the matrix diagram of $\pi$ has the form shown in Figure~\ref{fig2} above
where blank regions are empty, $A$ contains the lone entry $a$, $C$ contains the 
lone entry $c$, and, if $b$ is key, $B$ consists, left to right, of $s,s+1,\dots,b$, 
an interval of integers, for some $s\le b$, while if $b$ is not key, $B$ consists of 
$b+r,b+r+1,\dots,b+t,b+r-1,b+r-2,\dots,b$, for some $t\ge r$. Note that the central region $B$, colored yellow, is square-shaped and, since $A$ and $C$ each contain one entry, so are $X$ and $Z$.

The construction of $\si_1$ and $\si_2$ to define $\psi$ can now be described more explicitly:
$\si_2$ is obtained from $a\tau_2$ by replacing $c$ with the length of $a\tau_2$ and leaving all other 
entries intact, and $\si_1$ is obtained from $\tau_1 c$ as follows. 
\begin{itemize}
\item Case $b$ is key. Here, $B$ consists of $s,s+1,\dots,b$ with $2\le s\le b$. 
Replace $a$ with 1, subtract $s-1$ from each 
entry in $X$ and from the sole entry $c$ in $C$ (to account for $s-1$ missing entries $\{b\}\cup(\tau_2 \backslash \{c\})$\,), and subtract $s-2$ from each of the 
entries $s,s+1,\dots,b-1$ in $B$ (to account for $\tau_2 \backslash \{c\}$).
\item Case $b$ is not key. Here  $B$ 
consists of $b+r,b+r+1,\dots,b+t,b+r-1,b+r-2,\dots,b$. Replace $a$ with $r+1$, 
subtract $b-1$ from each entry in $X$ and from the sole entry $c$ in $C$,  subtract 
$b-2$ from each of the entries $b+r,b+r+1,\dots,b+t-1$ in $B$, and then append 
$r,r-1,\dots,1$. 
\end{itemize}
\begin{prop}
For  $1\le k <j <n$, the mapping $\psi$ is a bijection.
\end{prop}
\begin{proof}
In view of the preceding paragraph, to retrieve $\pi$ from $\psi(\pi)=(\si_1,\si_2)$, we need to recover $n,i,j,k,a,b,c$, to determine 
whether $b$ is key or not, and to find $r,s,t$ in the appropriate cases. Then we need to distinguish in $\si_1$ the entries increased by $s-1$ from those increased by $s-2$ when $b$ is key, and to distinguish the entries increased by $b-1$ from those increased by $b-2$ when $b$ is not key.
Clearly $n=\v\,\si_1\v+\v\,\si_2\v-1$ and $j=\v\,\si_1\v$. 
The index $i$ is the position of 1 in $\si_1$. 
Also,  $k$ = (position in $\si_2$ of its maximum entry) + $j-1$, and 
$a$ is the first entry of $\si_2$. The cases (1) $\pi_j=b$ is key and $\pi_{j-1}<b$, (2) $\pi_j$ is key and $\pi_{j-1}$ is a \rl, (3) $\pi_j$ is not key, can be distinguished according as (1) $\si_1$ ends with an ascent, (2) $\si_1$ ends with a descent to a number other than 1, (3) $\si_1$ ends with a descent to 1.
In all cases, the entries in $\pi$ after $b$ other than $c$ are evident in $\si_2$.
\begin{itemize}
\item Case $b$ is key. Here, $s$ is retrieved as the length of $\si_2$ and $i$ is the position of 1 in $\si_1$. Write $\si_1$ as $\si_{11}\si_{12}$ where $\si_{12}$ is the longest increasing right factor of $\si_1$. Then $b=s-1+(\textrm{length of }\si_{12})$. In $\si_1$, add $s-1$ to each entry of $\si_{11}$ and to the last entry of $\si_{12}$ and add $s-2$ to all other entries of $\si_{12}$ to obtain $\tau_1 c$.  
\item Case $b$ is not key. Here, $b$ is retrieved as the length of $\si_2$, $r$ is the length of the longest right factor of $\si_1$ of the form $r(r-1)\dots 1$, and $i$ is the position of $r+1$ in $\si_1$. Write $\si_1$ as $\si_1'r\,(r-1)\dots 1$ and write $\si_1'$ as $\si_{11}'\si_{12}'$ where $\si_{12}'$ is the longest increasing right factor of $\si_1'$. Then $t=r +\v\mskip2mu \si_{12}' \v -1$. In $\si_1'$,  add $b-1$ to each entry of $\si_{11}'$, add $b-2$ to all but the last entry of $\si_{12}'$, and append $b+t,b+r-1,b+r-2,\dots,b+1$ to get $\tau_1$. Add $b-1$ to the last entry of $\si_{12}'$ to get $c$.
\end{itemize}
  Thus $\psi$ is invertible.
\end{proof}

\section{Putting it all together} 
It is well known that 123-avoiders are counted by the Catalan numbers $C_n=\binom{2n}{n}-\binom{2n}{n-1}$ and so the number of start-small 123-avoiders of length $n$ is $C_n-C_{n-1}$ for $n\ge 1$, with \gf $C(x)-1-xC(x)=x^2 C(x)^3$ where $C(x)=\sum_{n\ge 0}C_nx^n$ is the \gf  for the Catalan numbers. 
Say a \textit{$k$-list} is a list of length $k$ and
the \emph{weight} of a start-small 123-avoider is one less than its length. 

From the bijection $\phi$, the number of \ss \avs of length $n$ is \begin{eqnarray}\label{wt}
& & \sum_{k=1}^{n-1}\textrm{\,\# $k$-lists of \ss 123-avoiders of total length $n-1+k$}  \nonumber  \\
& = & \textrm{\# lists (of arbitrary length) of \ss 123-avoiders of total weight $n-1$}.
\end{eqnarray}

Recall that if \a is a class (species) of combinatorial structures with $a_{n}$ structures of weight $n$ ($n\ge 1$),
and a \ \textit{compositional} \a-structure of weight $n$ is one obtained by taking a composition $(n_{1},n_{2}, \dots ,
n_{k})$ of $n$ and forming a $k$-list of \a-structures of respective weights $n_{1},n_{2}, \dots ,n_{k}$, then 
the counting sequence $(b_{n})_{n\ge 1}$ for compositional \a-structures by total weight has \gf $B(x):=\sum_{n\ge 1}b_{n}x^{n}$ given by
\[
1+B(x) =\frac{1}{1-A(x)},
\]
where $A(x)$ is the \gf for \a-structures (this is the INVERT transform). 
Applying it to the class of nonempty \ss 123-avoiders with weight measured as ``length minus 1'' so that $A(x)=x C(x)^3$,
we get 
\begin{equation}\label{B}
B(x)=\frac{1}{1-xC(x)^3} -1 \, .
\end{equation}

From (\ref{wt}) and (\ref{B}), the \gf $G(x)$ for nonempty \ss \avs (with $x$ marking length) is $G(x) = x B(x)$. 
Now let $(u_{n})_{n\ge 0}=(1,1,2,6, \dots)$ and $(v_{n})_{n\ge 0}=(1,0,1,4, \dots)$ denote the counting sequences 
for \avs and \ss \avs respectively. Clearly, $v_{n}=u_{n}-u_{n-1}$ for $n\ge 1$, (consider deletion of the 
first entry from a \av on $[n]$ that starts $n$). So the \gfs $F(x)=\sum_{n\ge 0}u_{n}x^{n}$ and 
$G(x)=\sum_{n\ge 1}v_{n}x^{n}$ are related by $F(x)=(1+G(x))/(1-x)$.
Thus
\[
F(x)=\frac{1+G(x)}{1-x} = \frac{1 + \frac{x}{1- x\,C(x)^{3}} -x}{1-x}
\]
which, after expansion, agrees with Kotesovec's formula.

Losonczy \cite{losonczy} has counted permutations that 
avoid 3421, 4312 and 4321 or equivalently (by reversal) both of the patterns treated here and 1234.

\section*{Acknowledgement} 
I thank one anonymous referee for their helpful comments and another for numerous suggestions that greatly improved the presentation.

\end{document}